\documentstyle[11pt,amscd,amsfonts]{amsart}
\newtheorem{thm}{Theorem}[section]
\newtheorem{lem}[thm]{Lemma}
\newtheorem{cor}[thm]{Corollary}
\newtheorem{prop}[thm]{Proposition}
\newtheorem{notitle}[thm]{ }

\theoremstyle{definition}

\def \smash {\wedge}
\def \colim {\mathop{\mathrm {colim}}}
\def \L {\mathop{\mathrm {L}}}

\def \C {\mathbb{C}}
\def \H {\mathbb{H}}
\def \Q {\mathbb{Q}}
\def \R {\mathbb{R}}
\def \T {\mathbb{T}}
\def \Z {\mathbb{Z}}

\def \CA {{\cal A}}

\def \CC {{\cal C}}

\def \sc {spin^c}

\def \f {\phi}

\def \lam {\lambda}
\def \Lam {\Lambda}

\def \s {\sigma}

\def \fs {\frak{s}}

\begin{document}

\baselineskip.525cm

\title[Cohomotopy invariants: I]{A stable cohomotopy refinement of \\Seiberg-Witten invariants: I}
\author[Stefan Bauer]{Stefan Bauer} \address{Fakult\"at f\"ur Mathematik,
Universit\"at Bielefeld \newline PF 100131, 33501 Bielefeld, Germany}
\email{bauer@@mathematik.uni-bielefeld.de}
\author[Mikio Furuta]{Mikio Furuta} \address{Graduate School of Mathematical Sciences, University
of Tokyo \newline Tokyo 153-8914, Japan}
\email{furuta@@ms.u-tokyo.ac.jp}
\begin{abstract}The monopole map defines an element in an equivariant stable
cohomotopy group refining the Seiberg-Witten invariant. Part I discusses
the definition of this stable homotopy invariant and its relation to
the integer valued Seiberg-Witten invariants. A gluing theorem for
these invariants, proved in part II, gives new results on diffeomorphism 
types of decomposable four-manifolds. 
\end{abstract}
\maketitle


\section{Introduction\label{Intro}}

\noindent
In the year 1994, Nathan Seiberg and Edward Witten
found new differential invariants \cite{Witten} 
of four-manifolds. The basic geometric
ingredient in their construction was a space well known to physicists:
The moduli space of monopoles. This space is defined for a closed Riemannian 
four-manifold $(X, g)$, which is $K$-theory 
oriented (or equivalently equipped with both an 
orientation and a $\sc$-structure). Under certain genericity assumptions,
this moduli space is a smooth, orientable finite dimensional manifold and comes with
a natural complex line bundle. The Seiberg-Witten invariant is the characteristic
number obtained from these data by specifying an orientation of the moduli space. 

In \cite{Furuta}, 
a new point of view was introduced to monopolized four-manifold theory: 
Instead of the moduli space of monopoles, consider the monopole map!
The monopole map $\mu=\mu_{g}$ is an ${\T}$-equivariant map between certain 
affine Hilbert spaces. The moduli space of monopoles is obtained as
the quotient of the zero-set of $\mu$ by the action of the
group ${\T}$ of complex numbers of unit length.

The present article carries this point of view a little 
further. In a finite dimensional analogue, the map $\mu$ as a proper
map between finite dimensional Hilbert spaces would extend to the one-point
compactifications, thus defining an element in some equivariant version of
the stable homotopy groups of spheres. This stable homotopy element
by the Pontrijagin-Thom construction would encode the moduli space as 
an element in some sort of equivariant framed bordism group.

With a grain of salt, this concept actually carries 
over to the infinite dimensional setting.
Salt is provided by a theorem of A. S. Schwarz \cite{Schwarz}, 
which seems to have fallen into oblivion during the past decades. 
Applying his ideas to the monopole map leads to the following theorem.

\begin{thm} The monopole map $\mu$ defines an element in an equivariant    
        stable cohomotopy group 
         \[\pi^b_{{\T},H}(Pic^0(X);{\mathrm ind}({\cal D})),\] 
which is independent of the chosen Riemannian metric.
            For $b> {\mathrm dim}(Pic^0(X))+1$, a homology 
           orientation determines a homomorphism 
            of this stable cohomotopy group to $\Z$, which maps
            $[\mu]$ to the integer valued Seiberg-Witten invariant. 
        \end{thm}

Here $Pic^0(X)$ denotes the Picard torus $H^1(X;{\R})/H^1(X;{\Z})$. The
Dirac operator associated to the chosen $\sc$-structure defines a
virtual complex index bundle 
${\mathrm ind}(\cal D)$ over the Picard torus, and 
$b=b_+(X)$ denotes the dimension of a maximal linear subspace $P$
of the second de Rham group of $X$, on which the cup product pairing is positive definite.
An orientation of $P\times Pic^0(X)$ is called a homology orientation. It determines
an orientation of the monopole moduli space.
The suffix $H$, finally, stands for a universe for the ${\T}$-action, that is
a Hilbert space with an orthogonal ${\T}$-action. Its finite dimensional
invariant linear subspaces provide the suspension coordinates in the
construction of equivariant cohomotopy groups. 

\noindent
The stable cohomotopy invariant as formulated in the theorem 
does not capture all the
features of Seiberg-Witten theory. In the case $b=1, Pic^0(X)=0$, for example,
the above groups are vanishing. The Seiberg-Witten invariants, however, are
nonvanishing in general, but depend in a well understood manner on the 
Riemannian metric and on an additional perturbation parameter. 

Indeed it is possible to recapture this phenomenon in
the stable homotopy setting at the expense of loosing convenient 
algebraic structure: One may consider $\mu$ up to equivariant
homotopy relative to the fixed point set. In case 
$b> {\mathrm dim}(Pic^0(X))+1$
the resulting relative homotopy classes are 
one-to-one with the homotopy classes which 
form the stable cohomotopy groups above.
In general the sets of such stable relative homotopy classes don't 
admit a natural group
structure.  The relative homotopy class of $\mu$
as an element in such a set turns out to be independent of 
the parameters used in the definition of the Seiberg-Witten number.
The reason is that the restriction of the monopole map to the ${\T}$-fixed point
set is quite special: It is affine linear.

So, how can a discrete invariant (as the Seiberg-Witten number) of any sort of 
homotopy class jump along a continuous path of parameters? To find an
answer, one has to
carefully distinguish between parameters: 
One parameter, the metric, will move
$\mu$ within its relative homotopy class. 
The second, problematic
parameter is used in the definition of 
a map from the set of relative homotopy classes
to the integers, which associates to the class of $\mu$ a number.
Seiberg-Witten theory defines such maps by performing an oriented 
count of orbits
in the preimage of generic fixed points in the target space of $\mu$.
However, there is no natural choice of a generic fixed point. 
The affine linear map $\mu^{{\T}}$ on the fixed point set 
maps onto a real hyperplane which divides 
the ${\T}$-fixed point set in the target space into two ``chambers''.
The oriented count has different results for generic points 
in different chambers.

Let's illustrate  this phenomenon in a characteristic example:
View the spinning globe as a two-sphere with an ${\T}$-action and choose
the north pole as a base point. As a target space, take  a one-sphere with
trivial action and choose two points on this one-sphere as ``poles'', the 
north pole again as base point. Based equivariant maps from the spinning globe
to the one-sphere are determined by their restriction to a latitude, 
which as an arc is a contractible space. So there is only the 
trivial homotopy class of equivariant such maps. 

In contrast, consider equivariant maps, which take north and south pole to
north and south pole, respectively. Such a map basically wraps a 
latitude $n+{1\over 2}$ times around the one-sphere. 
This wrapping number classifies the ``relative homotopy'' class of such a map.
Choosing a generic point in the one-sphere, the oriented count of 
preimages in a fixed latitude defines in a natural way a map of the set of
relative homotopy classes to the integers. This oriented count, 
however, depends
on the choice of the generic point. It changes by $\pm 1$, if 
the generic point is chosen in the ``other half'' of the one-sphere.

In the case of a spin-manifold, the monopole map is $Pin(2)$-equivariant.
It thus defines an element in a $Pin(2)$-equivariant cohomotopy group. 
One can recover the $10\over 8$-theorem of \cite{Furuta} in this setup:
In the case of vanishing
first Betti number, there is a fixed point homomorphism from this $Pin(2)$-equivariant
group to the $0$-th stable stem, i.e. to $\Z$.  It is shown in \cite{Birgit} 
that surjectivity of this fixed point homomorphism is equivalent to the 
existence of equivariant maps as used in the proof of \cite{Furuta}.

Donaldson's theorem on the 
diagonalizability of the quadratic form on manifolds with 
negative definite intersection form 
follows from number theory and from a well known result
on ${\T}$-equivariant maps between spheres.

\noindent
{\it Addendum:}
The stable homotopy invariant in the theorem was found independently by
both authors. For various reasons, publication of the corresponding
preprints \cite{PartII} and \cite{Furuta2} was unduly delayed. As a consequence, 
we finally opted for a joint publication.

\section{Fredholm maps and stable homotopy}

\noindent
A {\it Fredholm map} $f:H'\to H$ in this paper will be 
a compact perturbation  of a linear
Fredholm operator between separable Hilbert spaces. 
This means that
$f$ is of the form $f=l+c$, where $l$ is linear Fredholm and 
the continuous map $c$ maps bounded sets to subsets of compact sets.

\noindent
A theorem, which goes back to A.~S.~Schwarz \cite{Schwarz}, associates 
to certain such Fredholm maps stable homotopy classes of maps between finite 
dimensional spheres: 
Let $D'\subset H'$ be a disc with boundary $\partial D'$. Two continuous maps 
$f_i:\partial D'\to H\setminus \{0\}$
are called {\it compactly homotopic relative to $l$}, if there is a continuous and compact 
map $c:\partial D'\times [0,1]\to H$ with $f_i=l+c_i$ for $i\in\{0,1\}$ and 
\[
f_t(h')=l(h')+c_t(h')\not= 0
\]
for all $t\in [0,1]$ and $h'\in \partial D'$. 

\begin{thm}\label{Svarc}
           The compact homotopy classes 
           of continuous Fredholm maps relative to $l$ are in one-to-one correspondence 
           with elements of the stable homotopy group
           $\pi^{st}_{\mathop{\mathrm ind} l}(S^0)$ of the sphere.
\end{thm}

\noindent
This correspondence can be described as follows:
Any compact map $c$ on the bounded set $\partial D'$ can be uniformly approximated by maps
$c_n$ mapping to finite dimensional linear subspaces $V_n\subset H$ containing
$im(l)^\perp$. The correspondence associates to $l+c$ the maps $f_n/||f_n||$ with 
$f_n$ being the restriction of
$l+c_n$ to $l^{-1}(V_n)\cap\partial D'$.
A detailed proof of this theorem can be found in
\cite{Berger}, p.257f.

\noindent  
In this paper, this concept will be used with a few modifications:
Firstly, we will consider equivariant maps, which are furthermore
parametrized over some space.
Secondly, we will consider Fredholm maps which extend 
continuously to maps $$f^+:{H'}^+\to H^+$$ 
between the one-point completed Hilbert spheres. 
Equivalently, we suppose $f$ to satisfy a boundedness condition:
The preimages of bounded sets are bounded. 

The second condition is convenient and avoids 
clumsy notation: 
The construction of \cite{Schwarz} considers
maps between pairs $(D',\partial D')\to (H, H\setminus \{0\})$ of spaces.
The boundedness condition makes it possible to instead consider maps ${H'}^+\to H^+$
between spheres equipped with a natural base point.
As in finite dimensions,
both points of view are equivalent for maps which satisfy the
boundedness condition. This is due to the 
fact that starting from the pairs 
$({H'}^+, +)$ and  $(H^+, +)$ one gets the pairs 
$(D', \partial D')$ and $(H, H\setminus\{0\})$  via homotopy
equivalence and excision with intermediary pairs 
$({H'}^+, {H'}^+\setminus{\buildrel \circ\over D'})$ and $(H^+, H^+\setminus 
\{0\})$.
In the monopole setup the boundedness condition is satisfied. 

\noindent
First a short discussion of the boundedness condition:
In finite dimensions this condition
is equivalent to $f$ being proper, i.e. $f$ is closed and the
preimage of any point in the target space is compact.   
Here is a proof that in the setting of Fredholm maps 
the boundedness condition implies properness:

\begin{lem}
           Let $l: H'\to H$ be a continuous 
           linear Fredholm map between Hilbert spaces
           and let $c: H'\to H$ 
           be a compact map. Then the restriction of the map $f=l+c$
           to any closed and bounded subset $A'\subset H'$ is proper.\\
           In particular, if preimages of bounded sets
           in $H$ under the map $f$ are bounded, then $f$ is proper and 
           extends to a proper map $f^+:{H'}^+\mapsto {H}^+$ between the one 
           point completions.
        \end{lem}
\begin{pf}
           Let
           $\rho: H'\to \ker l$ denote the orthogonal projection. Then $f|_{A'}$ factors
           through an injective, closed and thus proper map 
           $A'\to H\times \overline{c(A')}\times \overline{\rho(A')}$, 
           $a'\mapsto (l(a'),c(a'), \rho(a'))$, a homeomorphism $(h,s,e)\mapsto(h+s,s,e)$
           and the projection
           $H\times \overline{c(A')}\times \overline{\rho(A')}\to H$ 
           which is proper as the two extra factors are compact.
        \par\noindent
         Now we invoke the boundedness condition: As the preimages of points
           in $H$ are bounded, they are compact by what was already shown. 
           Let $h\in H$ be in the closure of $f(A')$, with $A'$ closed in 
           $H'$. By the boundedness condition, $h$ is already in the closure 
           of $f({A'}_0)$, where ${A'}_0$ is a bounded closed subset of $A'$.
           From the first part of the proof it follows that $h$ is contained 
           in $f({A'}_0)\subset f(A')$. Thus $f$ is proper. But properness extends 
           to $f^+$: If, for a closed $A'\subset {H'}^+$, the closure
           $\overline{f^+(A')}$ contains the point at infinity, then 
           $f^+(A')\cap H=f(A'\cap H')$ is unbounded. Since $f$ is a compact perturbation
           of a continuous linear map, $A'\cap H'$ is unbounded and 
           thus contains the point at
           infinity in its closure. In particular, $f^+$ is  closed and thus proper.
        \end{pf}

\noindent
Indeed, proper Fredholm maps need not satisfy the boudedness condition. An 
example is easily provided: 
Let $H$ be a separable Hilbert space with orthonormal 
basis $e_n$. Let $\phi: H\to [0,1]$ be a continuous function taking 
values in the unit interval, which is supported in the ball in $H$ 
of radius $1\over 2$ and such that $\phi(0)=1$. Then the map defined by
$$f(x)= x+\sum_{n=1}^\infty (n-1)\phi(x-ne_n)e_n$$
is proper, but the preimage of the unit ball is not bounded.

\noindent 
We are now going to associate to a Fredholm map satisfying the boundedness
condition a stable homotopy class of maps between spheres. The next lemma
will provide the technical foundations. Let's start with fixing notation:\\  
Let $W\subset H$ be a finite dimensional linear subspace and let
$W'=l^{-1}(W)$ be its preimage under the linear Fredholm map $l$.
Let $S(W^\perp)$ denote the unit sphere in the orthogonal complement 
$W^\perp$ of $W$. 
As in finite dimensions, the inclusion ${W}^+\to {H}^+\setminus S(W^\perp)$
is a deformation retract. The retracting map $\rho_W$ can be described as
follows: The one-point completed Hilbert space $H^+$ identifies 
with the unit sphere $S({\R}\oplus H)=S({\R}\oplus W\oplus W^\perp)$ 
in ${\R}\oplus H$ via the map $h\mapsto (|h|^2+1)^{-1}(|h|^2-1, 2h)$. 
In this identification, the subspace $W^+$ maps to the "equatorial" subsphere 
$S({\R}\oplus W\oplus 0)\subset S({\R}\oplus W\oplus W^\perp)$ and $S(W^\perp)$ maps to
the complementary "polar" subsphere $S(0\oplus 0\oplus W^\perp)$. The retracting
homotopy shrinks the latitutes in $S({\R}\oplus W\oplus W^\perp)\setminus
S(0\oplus 0\oplus W^\perp)$ to the equator. The retraction $\rho_W$ has the 
following property: For $h\in H\setminus W^\perp$, the vector $\rho_W(h)$ differs
from the orthogonal projection $pr_W(h)$ to $W$ by a positive scalar factor
$\rho_W(h)= \lambda (h)pr_W(h)$.

\begin{lem}\label{Einhaengung} There are finite dimensional linear subspaces
$V\subset H$, such that the following statements hold:
 \begin{enumerate}
          \item 
             The subspace $V$ spans, together with the image $Im (l)$ of the 
             linear Fredholm map $l$, the Hilbert space $H=Im(l)+V$.
          \item
             For $W\supset V$ with $W= U\perp V$, 
             the restricted map $f|_{{W'}^+}:{W'}^+\to H^+$
             misses the unit sphere $S(W^\perp)$ in the orthogonal
             complement of $W$.
          \item The maps $\rho_{W}{f|_{{W'}^+}}$ and 
             ${\mathrm id}_{U^+}\smash \rho_{V}f|_{{V'}^+}$ 
             are homotopic as pointed maps
             \[{W'}^+\cong{U}^+\smash {V'}^+\to {U}^+\smash V^+= W^+.
             \]
       \end{enumerate}   
Indeed, if $H$ is separable, then the subspaces $V$ satisfying these three conditions
are cofinal in the direct system of finite dimensional subspaces in $H$. 
\end{lem}

\begin{pf}    The preimage $f^{-1}(D)$ of the unit disc $D$ in $H$
              is bounded in $H'$. So the closure $C$ of its image under the
              compact map $c$ is compact in $H$. Cover $C$ by finitely many balls with 
              radius $\varepsilon\le {1\over 4}$, 
              centered at points $v_i$. Together with the orthogonal complement 
              to the image of the linear Fredholm map $l$, these points $v_i$
              span  a finite dimensional linear subspace $V$ of $H$. 
              Let's check the second condition: Suppose $w\in S( W^\perp)$ is
              in the image of ${f|_{{W'}^+}}$. Then 
              $f^{-1}(w)\cap {W'}^+\subset f^{-1}(D_1(H))$
              will be mapped by ${f|_{W'}}=(l+c)|_{W'}$ to a subspace of $W+C$.
              So $w$ will be contained both in $S(W^\perp)$ and $W+C$. However, these
              two subsets of $H$ are at least $1-\varepsilon\ge {3\over 4}$ apart.
             \par\noindent
              We will identify $W'$ with the orthogonal sum $U\perp V'$ via the
              map $$w'\mapsto (l\circ(1-pr_{V'})(w'),pr_{V'}(w')).$$
             \par\noindent
              To prove the last claim, it suffices to show that 
              ${\mathrm id}_{{U}^+}\smash \rho_{V}f|_{V'}^+$ and
              $f|_{W'}^+$ are homotopic as maps ${W'}^+\to H^+\setminus S(W^\perp)$.
              Let $D'\subset H'$ be a disk, centred at the origin,
              which contains the preimage $f^{-1}(D)$ 
              of the unit disk in $H$. Consider the homotopy 
              $h:D'\times [0,3]\to H^+\setminus S(W^\perp),$
              defined by:
              $$
                h_t=\left\{
                          \begin{matrix}                          
                l+((1-t)id_H+t\cdot pr_V)\circ c   &   \text{for $0\leq t\leq 1$,}\\
                l+ pr_V\circ c\circ((2-t)id_{V'}+(t-1)pr_{V'}  
                                              &   \text{for $1\leq t\leq 2$,}\\
                pr_{U}\circ l+((3-t)pr_V+(t-2)\rho_V)\circ (l+c)\circ pr_{V'}   
                                              &   \text{for $2\leq t\leq 3$.}\\
                           \end{matrix}
               \right. 
              $$
              Note that the image during the homotopy stays within an 
              $\varepsilon$-neighbourhood of $W$. The homotopy is chosen 
              in such a way that the image of the sphere 
              $S'=\partial D'\cap W'$ during the 
              homotopy stays away not only from the unit sphere $S(W^\perp)$ in
              $W^\perp$, but from the whole of $W^\perp$. Before we check this,
              let's consider the consequences: Since $H^+\setminus (D\cap W^\perp)$
              is contractible, the homotopy $h_t$ can be extended to the complement
              of $D'\cap W'$ in ${W'}^+$, thus defining a homotopy as claimed.
             \par\noindent
              Let $s'$ be an element in the sphere $S'$. We will 
              track the path its image will take during the homotopy. At starting
              time, it is mapped to $f(s')$, which is of norm greater or 
              equal to $1$, and furthermore, in an $\varepsilon$-neighbourhood
              of $W$.  In particular, its distance form $W^\perp$ is at least
              $1-\varepsilon\geq \frac34$. During the first part of the homotopy,
              the image will move at most a distance of $\varepsilon$, so it 
              will definitely stay away from $W^\perp$.
             \par\noindent
              From time $t=1$ on one has $pr_{U}\circ h_t(s')=pr_{U}\circ l(s')$.
              Since $pr_{U}(W^\perp)=0$ by definition, we are reduced to checking
              the case $pr_{W'}\circ l(s')=0$, that is for $s'\in S'\cap V'$.
              But for such an element, the image during the second part of the
              homotopy stays fixed and during the third part moves on a straight
              line between $pr_V(f(s')) $ and $\rho_V(f(s'))$, which are 
              nonzero vectors in $V$, differing by a positive real factor. This
              concludes the proof of \ref{Einhaengung}.  
\end{pf}

\noindent
In particular, the restrictions $f|_{l^{-1}(V)}$ to finite dimensional linear 
subspaces $V\subset H$ as in \ref{Einhaengung} together define an element in the
colimit of pointed homotopy classes
$$
[f]=\colim_{V\subset H}[(f|_{l^{-1}(V)})^+]\in 
\colim_{V\subset H}[(l^{-1}(V))^+,H^+\setminus S(V^\perp)].
$$
The homotopy equivalences $V^+\subset (H^+\setminus S(V^\perp))$
combine to an isomorphism
$$
\pi^{st}_{\mathop{\mathrm ind} l}(S^0)=\colim_{V\subset H}\ [(l^{-1}(V))^+,V^+]
             \ {\mathop{\rightarrow}^\sim}\ \ 
             \colim_{V\subset H}\ [(l^{-1}(V))^+,H^+\setminus S(V^\perp)].
$$
In this way $[f]$ can be identified as an element in the stable homotopy group
$\pi^{st}_{\mathop{\mathrm ind} l}(S^0) $:

\begin{cor}
Let $f=l+c:H'\to H$ be a compact perturbation of the linear Fredholm map $l$
such that the preimages of bounded sets under the map $f$ are bounded. Then
$f$ defines an element $[f]\in  \pi^{st}_{\mathop{\mathrm ind} l}(S^0) $.
\qed\end{cor}
\noindent

\noindent
In the construction above the linear map $l$ seems to play an essential r\^ole.
In fact it will turn out that the homotopy class $[f]$ basically is independent of
the choice of decomposition of $f$ as a sum $f=l+c$. In order to show this, we 
will have to consider a parametrized version of the above situation and reach 
back some way:

\noindent
Let $Y$ be a finite CW-complex. The group $KO^0(Y)$ can be described
as follows (cf. \cite {Segal}):

\noindent
A (real) Hilbert bundle over $Y$ is a locally trivial fiber bundle with fiber
a separable Hilbert space $H$, whose structure group is the group of linear 
isometric bijections, equipped with the norm topology.
A {\it cocycle} $\lambda=(E',l,E)$ over $Y$ consists of two Hilbert bundles
over $Y$ and a Fredholm morphism $l:E'\to E$ between them. Here a Fredholm 
morphism is a continuous map which is fiber preserving and fiberwise linear
Fredholm over $Y$. Two cocycles $\lambda_i$ over $Y$ for $i\in\{0,1\}$ are 
homotopic, if there is a cocycle $\lambda$ over $Y\times [0,1]$ such that the restriction
$\lambda|_{Y\times\{i\}}$ is isomorphic to $\lambda_i$. A cocycle $(E',l,E)$ is
trivial, if $l$ is invertible. Two cocycles $\lambda_0$ and $\lambda_1$ are 
equivalent, if there is a trivial cocycle $\tau$ such that $\lambda_0\oplus \tau$
and $\lambda_1\oplus \tau$ are homotopic. The group $KO^0(Y)$ is the set of 
equivalence classes of cocycles with addition given by the Whitney sum of cocycles.
\par\noindent
Let $f:E'\to E$ be a continuous map between Hilbert bundles of the form
$f=l+c$, where $\lambda=(E',l,E)$ is a cocycle over $Y$ and $c$ is fiber
preserving and compact, i.e. maps bounded disk bundles in $E'$ to subspaces in $E$,
which are proper over $Y$. Let's call such a map $f$ a {\it Fredholm map 
over} $Y$. The boundedness condition in this parametrized situation reads:
The preimages of bounded disk bundles are contained in bounded disk bundles. An
equivalent condition is: The Fredholm map over $Y$ extends to the fiberwise one-point
completions of $E'$ and $E$.
\par\noindent
Every Hilbert bundle over the compact space $Y$ is trivial, i.e. $E\cong Y\times H$
by the theorem of Kuiper \cite{Kuiper}.  The boundedness condition on $f$
thus translates to the condition that the composed map 
$pr_H\circ f: E'\to H$ extends to the one-point
completions, defining a continuous map
$$
(pr_H\circ f)^+:T(E')\to H^+
$$
from the Thom space of the Hilbert bundle $E'$ to the Hilbert sphere $H^+$.   
\par\noindent 
The stage is now set for the definition of stable cohomotopy groups with coefficients:  
Let $\lambda$ be a finite dimensional virtual vector bundle over $Y$. Suppose we
are given a presentation $\lambda=F_0-F_1$ with vector bundles $F_i$ such that
$F_1\cong Y\times V$ is a trivial vector bundle with $V$ a finite dimensional
linear subspace of a Hilbert space $H$. 
With $TF_0$ denoting the Thom space of the bundle $F_0$,  
stable cohomotopy groups are defined as the colimits 
\begin{eqnarray*}
    \pi_H^n(Y;\lambda)&=&
       \colim_{U\subset V^{\perp}}\ [{U}^+\smash TF_0\ ,\ {U}^+\smash V^+\smash S^n]\\
           &=&\colim_{W\subset H}\ [W^+\smash T\lambda\ ,\ W^+\smash S^n]
\end{eqnarray*}
of pointed homotopy classes of maps, where the colimits are over the finite dimensional 
linear subspaces $U\subset V^{\perp}\subset H$ and $W=U+V\subset H$, 
respectively. Here the connecting
morphism for $W\subset W_1$  with $U_1=W^\perp\cap W_1$ is the suspension map
(${\mathrm id}_{{U_1}^+}\smash \, .$). 
The symbol $T\lambda$ stands not anymore for a space, but for a spectrum.

\noindent
The reason for keeping the Hilbert space $H$ in the notation lies in the equivariant
version:
For a compact Lie group $G$ we fix 
a $G$-universe $H$, i.e. a real Hilbert space $H$ equipped with an orthogonal
$G$-action such that $H$ contains the trivial representation and, furthermore,  
the space of equivariant morphisms $Hom_G(V,H)$ for a real $G$-module $V$
either is zero or infinite 
dimensional.   Let $\lam=F_0-F_1$ be a virtual equivariant vector bundle over a 
finite $G$-CW complex $Y$ such that $F_1\cong Y\times V$ is a trivial bundle 
with $V\subset H$ a finite dimensional $G$-subrepresentation. 
Stable equivariant cohomotopy groups are the colimits 
\begin{eqnarray*}
    \pi_{G,H}^n(Y;\lam)&=&
        \colim_{U\subset V^{\perp}}\ [{U}^+\smash TF_0\ ,\ {U}^+\smash V^+\smash S^n]^G\\
           &=&\colim_{W\subset H}\ [W^+\smash T\lam\ ,\ W^+\smash S^n]^G
\end{eqnarray*}
of pointed equivariant homotopy classes of maps, where the colimit now is 
over the finite dimensional 
subrepresentations $U\subset V^{\perp}\subset H$ and $W=U+V\subset H$, 
respectively.
This definition of stable equivariant cohomotopy groups differs a little
from the usual one as we allow for coefficients $\lam$ in the equivariant
$KO$-group $KO^0_G(Y)$ and our universe 
$H$ need not contain all irreducible representations. 

\noindent
Let $f:E'\to E$ be a G-equivariant Fredholm map between G-Hilbert space bundles
over the finite $G$-CW complex $Y$ such that $E\cong Y\times H$ is a trivialised bundle.
Let $f=l+c$ be a presentation of $f$ as a sum of a linear Fredholm morphism and
a compact map.   For sufficiently large linear $G$-subspaces 
$V\subset H$, the cocycle $\lam=(E',l,E)$ admits a presentation as virtual index bundle   
\[\lam=F_0(V)-F_1(V) \] with equivariant vector bundles $F_0(V) =
({\mathrm pr}_H\circ l)^{-1}(V)\subset E'$ and $F_1(V)=Y\times V$. 
The following lemma parallels \ref{Einhaengung}. Its proof is omitted, as it is 
almost verbatim the same.

\begin{lem} There exist finite dimensional linear $G$-subspaces 
             $V\subset H$ such that the following hold:
      \begin{enumerate}
          \item 
             For every $y\in Y$, the subspace $V$ is mapped onto
             ${\mathrm coker}\ (l_y:E_y'\to H)$. In particular, 
             $F_0(V)$ is a bundle over $Y$ and $\lam=F_0(V)-F_1(V)$
             represents the virtual index bundle ${\mathrm ind}(l)$.
          \item 
             For any $G$-linear $W=W'+V$ with $W'\subset V^\perp$, the restricted map
             $f(W)^+=(pr_H\circ f)|_{F_0(W)}^+:TF_0(W)\to H^+$ 
             misses the unit sphere $S(W^\perp)$.  
          \item The maps $\rho_{W}f(W)^+$ and 
             ${\mathrm id}_{{W'}^+}\smash \rho_{V}f(V)^+$ are $G$-homotopic 
             as pointed maps
             \[F_0(W)^+\cong{W'}^+\smash F_0(V)^+\to {W'}^+\smash V^+=W^+.
             \qed\]
       \end{enumerate}   
 \end{lem}

\begin{thm}\label{Freddy}
  An equivariant Fredholm map $f=l+c:E'\to E$ between
             G-Hilbert space bundles over $Y$ with $E\cong Y\times H$, which extends
             continuously to the fiberwise one-point completions,
             defines a stable cohomotopy Euler class
             \[ [f]\in\pi^0_{G,H}(Y;\mathop{\mathrm ind} l).\]
             This Euler class is independent of the presentation of $f$ as a sum.  
                 
\end{thm}

\begin{pf}
The only statement left to prove is the final one. Note that the restriction maps
$$\pi^n_{G,H}(Y\times [0,1],\lam)\to  \pi^n_{G,H}(Y\times \{i\},\lam|_{Y\times\{i\}}))$$
are isomorphisms. Thus a homotopy of  cocycles naturally 
induces an isomorphism
of the corresponding cohomotopy groups. (An extension of this statement to
equivalences of cocycles needs further discussion of universes; it seems
unnecessary in the present context.) 
If $f=l+c=l'+c'$ are two different presentations
as a sum, then the constant homotopy $f=f_t=(1-t)(l+c)+t(l'+c')$ 
defines an Euler class in the cohomotopy group of $Y\times [0,1]$, which restricts
for $i\in\{0,1\}$ to the Euler classes defined via the respective presentations
of $f$.
\end{pf}

\begin{notitle} {\bf Remarks.}\end{notitle} 
\begin{itemize}
              \item 
                Indeed any element in $\pi^0_{G,H}(Y;\mathop{\mathrm ind} l)$ can be realized by
                a map between Hilbert space bundles satisfying the boundedness condition: 
                Take a finite 
                dimensional representative $p^+:TF_0^+\to V^+$. After
                possibly stabilizing further, this map is homotopic to one where the preimage 
                of the basepoint consists only of the base point. Now take the smash product 
                with the identity on an infinite dimensional Hilbert sphere and
remove the base point. 
              \item
                The stable cohomotopy Euler class has been defined and investigated by Crabb and Knapp
                \cite{CrabbKnapp}. It is related to the standard Euler class the following way:
                A section of an oriented vector bundle $\xi$ over $Y$ 
                can be regarded as a map $\sigma:Y\times {\R}^0\to\xi$. Choosing an bundle isomorphism
                $\xi\oplus \eta\cong Y\times {\R}^n$, this section and the projection to
                the fibers of a trivialized bundle together define a map
                $(\sigma+id_\eta)^+:\eta^+\to(Y\times{\R}^n)^+\to S^n$ and thus an element of 
                $\pi^0(Y;-\xi)$. The choice of a Thom class $u\in H^r(Y;\xi)=H^r(D\xi, S\xi)$
                corresponds to choosing an orientation of $\xi$. The standard Euler
                class is defined by $e(\xi)=\sigma^*(u)\in H^r(Y)$. A generator
                $1\in\tilde{H}^0(S^0)$ gives rise to 
                the Hurewicz map $\pi^0(Y;-\xi)\to H^0(Y;-\xi)$,  which associates to a 
                stable pointed map $\sigma:T(-\xi)\to S^0$ the image $\sigma^*(1)$.
                Using the cup product pairing $H^*(Y;-\xi)\times H^*(Y;\xi)\to H^*(Y)$,
                the singular cohomology Euler class and the stable cohomotopy one
                are related by       
                \[ e(\xi)=\sigma^*(1)\cdot u.\] 
\item
The approach of \ref{Svarc} and the one outlined above obviously are closely related: 
If $f=l+k:H'\to H$ admits {\it a priori}
estimates and $D'\setminus\partial D'\subset H'$ contains $f^{-1}(0)$, then its compact 
homotopy class in ${\cal C}^0_{l}(\partial D', H\setminus\{0\})$ corresponds to $[f]\in 
\pi^0_{H}(pt;\mathop{\mathrm ind} l)=\pi^{st}_{{\mathrm ind } l}(S^0)$.

\end{itemize}


\section{The monopole map}

\noindent
Let $S^+$ and $S^-$ denote the Hermitian rank-2 bundles associated to the given
$Spin^c$ structure on $X$ and let $L$ denote their determinant line bundle. 
Clifford multiplication \hbox{$T^*X\times S^\pm\to S^\mp$} defines a linear map 
$\rho:\Lam^2\to \mathop{\mathrm{End}}_{\C}(S^+)$ from the bundle of 2-forms
to the endomorphism bundle of the positive spinor bundle. The kernel of this
homomorphism is the subbundle $\Lam^-$ of anti-selfdual 2-forms. Its image 
is the subbundle of trace-free Hermitian endomorphisms.
For a $\sc$-connection $A$, denote by $D_A:\Gamma(S^+)\to\Gamma(S^-)$ 
its associated Dirac operator. The monopole map $\tilde{\mu}$ is defined 
for triples $(A,\f,a)$ of a $\sc$-connection $A$, a positive spinor $\f$ 
and a 1-form $a$ on $X$ by
\begin{eqnarray*}            
            \widetilde{\mu}:{\CC}onn\times\left(\Gamma(S^+)\oplus \Omega^1(X)\right)&\to&
            {\CC}onn\times\left(\Gamma(S^-)\oplus \Omega^+(X)\oplus H^1(X;{\R})\oplus
             \Omega^0(X)/{\R}\right)\\
            (A,\phi, a)&\mapsto&(A,\,D_{A+a}\phi,\,F_{A+a}^+-\s(\phi),\,a_{harm},
            \,{\mathrm d}^*a).
       \end{eqnarray*}
Here $\s(\phi)$ denotes the trace free endomorphism $\phi\otimes\phi^*-
{1\over 2}\vert\phi\vert^2\cdot \mathop{\mathrm{id}}$ of $S^+$, considered via the
map $\rho$ as a selfdual 2-form on $X$.
As a map over the space ${\CC}onn$ of $\sc$-connections, the monopole map is
equivariant with respect to the action of the gauge group ${\cal G}=map(X,{\T})$.
This group acts on spinors via multiplication with $u:X\to {\T}$, on connections
via addition of $iu{\mathrm d}u^{-1}$ and trivially on forms. 
Fixing a base point $\ast\in X$,
the based gauge group ${\cal G}_0$ is obtained as the kernel of the evaluation
homorphism $map(X,{\T})\to {\T}$ at $\ast$. 
\par\noindent
Let $A$ be a fixed connection. The subspace $A+ ker({\mathrm d})\subset {\CC}onn$
is invariant under the free action of the based gauge group with quotient space
isomorphic to \[Pic^0(X)= H^1(X;{\R})/H^1(X;{\Z}).\] 
Let $\CA$ and $\CC$ denote the
quotients
\begin{eqnarray*}
 {\CA}&=&(A+ker\ {\mathrm d})\times\left(\Gamma(S^+)\oplus \Omega^1(X)\right)/{\cal G}_0\\
{\CC}&=&(A+ker\ {\mathrm d})\times
                \left(\Gamma(S^-)\oplus \Omega^+(X)\oplus H^1(X;{\R})\oplus
                \Omega^0(X)/{\R}\right)/{\cal G}_0
\end{eqnarray*}
by the pointed gauge group. Both spaces are bundles over $Pic^0(X)$ and the quotient
\[\mu=\widetilde{\mu}/{\cal G}_0:{\CA}\to{\CC}\]
of the monopole map is a fiber preserving, ${\T}$-equivariant map over $Pic^0(X)$.

\noindent
For a fixed $k > 4$, consider the fiberwise ${\L}^2_k$ Sobolev completion ${\CA}_k$
and the fiberwise ${\L}^2_{k-1}$ Sobolev completion ${\CC}_{k-1}$
of $\CA$ and $\CC$.
The monopole map extends to a continuous map 
$\mu=\mu_k:{\CA}_k\to {\CC}_{k-1}$ over $Pic^0(X)$.  
It is the sum $\mu=l+c$ of the linear Fredholm 
map $l=D_A \oplus d^+\oplus \mathop{\mathrm pr}_{harm}\oplus {\mathrm d}^*$ 
and a term $c:(\phi,a)  \mapsto (0,F_A^+,0,0)+(a\cdot\phi, -\s(\f),0,0)$.
This map $c$ is compact as the sum of the constant map $F^+_A$ and the 
composition of a multiplication map 
${\CA}_k\times{\CA}_k\to{\CC}_k$, which is continuous for $k>2$, and 
a compact restriction map ${\CC}_k\to{\CC}_{k-1}$. 
The following statement and its proof
are  only slight variations of standard ones
in Seiberg-Witten theory, compare e.g. \cite{KroMro}:

\begin{prop}\label{eigentlich}
           Preimages $\mu^{-1}(B)\subset {\CA}_k$ of bounded disk bundles 
           $B\subset{\CC}_{k-1}$ are contained in bounded disk bundles.
        \end{prop}

\begin{pf} 
It is sufficient to prove this fiberwise for the Sobolev completions
of the restriction of the monopole map to the space $\{A\}\times (\Gamma(S^+)
\oplus ker({\mathrm d}^*))$, which maps to $\{A\}\times (\Gamma(S^-)\oplus
\Omega^+(X)\oplus H^1(X;{\R}))$. (Note that the monopole map is a linear isomorphism
on the additional factor; in the stable homotopy picture the restricted map and the 
full monopole map thus differ only by infinite dimensional, but otherwise 
irrelevant suspension.)
Using the elliptic operator $D=D_A+d^+$ and its adjoint, define the 
$L^2_k$-norm via the scalar product on 
           the respective function spaces through
          \begin{eqnarray*}
                (\, .\, ,\, . \,)_i=  
                (\, .\, ,\, . \,)_0   +   (D\, .\, ,D\, .\,)_{i-1},\,\,\,\,\,\,\,\,\,\,\,\,\,\,
                (\, .\, ,\, . \,)_0    =   (\, .\, ,\, . \,)    =    \int_X<\, .\, ,\, . \,>*1
          \end{eqnarray*}
          The norms  for the $L^p_k$-spaces are defined correspondingly.
          Let $ \mu(A,\phi,a)   =   (A,\varphi , b , a_{harm} ) \in {\CC}_{k-1}  $ 
           be bounded by some constant $R$. The Weitzenb\"ock formula for the Dirac operator  
           associated to the connection $A+a=A'$ reads
           \[
                D^*_{A'} D_{A'}  =  \nabla^*_{A'} \nabla_{A'}  +  {1\over4}s  -  {1\over2}F_{A'}^+,
           \]
           with $s$ denoting the scalar curvature of $X$.
           As a consequence, there is a pointwise estimate:
            \begin{eqnarray*}
                    \Delta|\phi|^2&=&
                    2<\nabla^*_{A'}\nabla_{A'}\phi,\phi>-2<\nabla_{A'}\phi,\nabla_{A'}\phi>\\
                    &\le&2<\nabla^*_{A'}\nabla_{A'}\phi,\phi>\\
                   &=& 2<D^*_{A'}D_{A'}\phi-{s\over4}\phi+{1\over2}F_{A'}^+\phi,\phi>\\
                   &=& <2D_{A'}^*\varphi-{s\over2}\phi+(b+\sigma(\f))\phi,\phi>
              \end{eqnarray*}               
             In particular, there are inequalities
            \begin{eqnarray*}
                    \Delta|\phi|^2+{s\over2}|\phi|^2+{1\over2}|\phi|^4&\le&
                    <2D_A^*\varphi,\phi>+<2a\cdot\varphi,\phi>+<b\phi,\phi>\\
                     &\le& 2(||D^*_A\varphi||_{{\L}^\infty} + ||a||_{{\L}^\infty}         
                     ||\varphi||_{{\L}^\infty})\cdot |\phi| + ||b||_{{\L}^\infty}\cdot |\phi|^2\\
                     &\le&c_1\Big((1+||a||_{{\L}^\infty})||\varphi||_{{\L}^2_{k-1}}\cdot|\phi|+
                     ||b||_{{\L}^2_{k-1}}\cdot|\phi|^2\Big),
             \end{eqnarray*}       
            with a constant $c_1$ by applying the Sobolev embedding theorems. To get a bound for
            the remaining term, use a Sobolev estimate 
            $||a||_{{\L}^\infty}\le c_2||a||_{{\L}^p_1}$ for some $p>4$ and the elliptic
            estimate $||a||_{{\L}_1^p}\le c_3(||d^+a||_{{\L}_0^p}+||a_{harm}||)$. 
            Combination with
             the equality $d^+a=b-F^+_A+\sigma(\f)$ then leads to an estimate
             \begin{eqnarray*}
                ||a||_{{\L}^\infty}&\le& c_4(||a_{harm}||+||b||_{{\L}^p_0}+||F^+_A||_{{\L}^p_0}+
                                     ||\sigma(\f)||_{{\L}^p_0})\\
                                  &\le& c_5(||a_{harm}||+||b||_{{\L}^2_{k-1}}+
                                              ||F^+_A||_{{\L}^p_0}+
                                              ||\f||^2_{{\L}^\infty}).
             \end{eqnarray*}
             At the maximum of $|\f|^2$, its Laplacian is non-negative. So, putting  
             everything together,
             one obtains a polynomial estimate of the form
             $$
               ||\f||_{{\L}^\infty}^4 \le cR\Big((1+R)||\f||_{{\L}^\infty}+
                                          ||\f||^2_{{\L}^\infty}+||\f||^3_{{\L}^\infty}\Big)+
                                          ||s||_{{\L}^\infty}||\f||^2_{{\L}^\infty}.
             $$  
           \noindent
            Combining the last two estimates, one obtain bounds 
            for the ${{\L}^\infty}$-norm and         
           a fortiori for the ${\L}^p_0$-norm of $(\f,a)$ for every $p\ge 1$. \par
         \noindent
           Now comes bootstrapping: For $ i \le k$, assume inductively 
           ${\L}_{i-1}^{2p}$-bounds on $ ( \phi , a )$ with $p=2^{k-i}$.  
           To obtain ${\L}_{i}^p$-bounds, compute:
           \begin{eqnarray*}     ||(\phi,a)||^p_{{\L}^p_i}  -   ||(\phi,a)||^p_{{\L}^p_0}
                            & =& 
                                 ||(D_A \phi,d^+a)||^p_{{\L}^p_{i-1}}                       \\
                            &=& 
                                 ||(\varphi,b,a_{harm})||^p_{{\L}^p_{i-1}}   +      
                                 ||(a \phi,-F_A^+-\s(\f))||^p_{{\L}^p_{i-1}}.
           \end{eqnarray*}
           The latter equality holds as $ D_{A'} = D_A + a$. 
           The summands in the last expression are bounded by the assumed   
           ${\L}_{i-1}^{2p}$-bounds on $ ( \phi , a )$.
\end{pf}

\noindent
The proposition in particular implies that 
the assumptions of \ref{Freddy} are satisfied for the monopole
map $\mu$. The conclusion is spelled out in the following

\begin{cor}The monopole map defines an element $[\mu]$
in the stable cohomotopy group 
\[ \pi^0_{{\T}, H}(Pic^0(X);\lam)=\pi^b_{{\T},H}(Pic^0(X);{ind}(D)),\]
where $H$ is a Sobolev completion of 
the sum  $\Gamma(S^-\oplus \Lam^2_+(T^*X))$ 
of the vector spaces of negative spinors and selfdual two-forms on $X$.
The virtual
index bundle $\lam=ind(D)\ominus H_+$ is the difference of the complex
virtual index bundle of the Dirac operator over $Pic^0(X)$ and the trivial
bundle  $H_+$ with fiber $H^2_+(X;{\R})$, which for a chosen metric on $X$
may be viewed as the space
of selfdual harmonic two-forms. The ${\T}$-action on $ind(D)$ is 
by multiplication with complex numbers and on $H_+$ is trivial.
\qed
\end{cor}

\noindent
There is a comparison map from the stable equivariant cohomotopy
group above to the integers, which relates the element defined
by the monopole map with the integer valued
Seiberg-Witten invariant associated to it:

\begin{prop}\label{standardSW} Let $X$ be a closed 4-manifold with 
             $b=b_+>b_1+1$. The choice of a homology 
             orientation (i.e. an orientation of 
             $H^1(X;{\R})\oplus H^2_+(X;{\R})$)
             then determines a homomorphism 
             $t:\pi_{{\T},{H}}^b(Pic^0(X); ind(D))\to \Z$,
             which maps the class of the monopole map to the integer valued
             Seiberg-Witten invariant.
\end{prop}

\begin{pf} Any element in $\pi_{{\T},{H}}^b(Pic^0(X); ind({\cal D}))$ 
           is represented by a pointed equivariant map 
           \[\mu: TF\to V^+\]
           from the Thom space of a bundle $F$ over $Pic^0(X)$ to a sphere 
           $V^+=(V'\oplus H^2_+(X;{\R}))^+$, where 
           $F-Pic^0(X)\times V'$ represents the equivariant virtual index 
           bundle of the Dirac operator. The integer Seiberg-Witten invariant is
           constructed as follows. After possible perturbation of the map $\mu$,
           the ${\T}$-fixed point set $TF^{{\T}}$ is mapped to  a subspace of $V^{{\T}+}$ of 
           codimension at least  $b-b_1\geq 2$. After perturbing further, the preimage of a 
           generic point in the complement is a manifold $M$ with a free ${\T}$-action.
           The homology orientation, together with a standard orientation of complex vector
           spaces, defines a relative orientation for the pair $F, V$ and thus an orientation
           of $M$ and on the manifold $M/{{\T}}$. The dimension of $M$ is $ind_{\R}(D)+b_1-b_2$.
           Now suppose $M$ is a manifold of odd dimension $2k+1$ (otherwise the SW-number is zero). 
           Then the 
           Seiberg-Witten number is the evaluation of the Euler class of the complex vector bundle
           $(M\times{\C}^k)/{{\T}}$ over $M/{{\T}}$ at the fundamental class. 
           
           \noindent Equivalently, one could start with 
           the map $\gamma({\C})^k\circ \mu$, where $\gamma({\C})^k:{V}^{+}\to (V\oplus {\C}^{k})^+$ is the 
           one-point compactified inclusion of vector spaces. After perturbing this map 
           equivariantly as before, the preimage of a generic fixed point would be 
           a finite number of oriented ${\T}$-orbits. The oriented count of it is the Seiberg-Witten
           number again. This is basically the same construction as above, as the Euler class of
           a vector bundle over a manifold is the Poincar\'e dual of a generic section.
           
           \noindent
           Here is another equivalent description of this map in purely algebraic topological
           terms. In the long exact sequence of stable homotopy groups of equivariant maps
           associated to the pair $(TF, TF^{{\T}})$
           $$\begin{CD} 
                      \{\Sigma TF^{{\T}}, V^+\}_{{\T},H} \to \{(TF,TF^{{\T}}), (V^+,\emptyset^+)\}_{{\T},H} 
                        \to \{TF,V^+\}_{{\T},H} \to \{TF^{{\T}},V^+\}_{{\T},H}
           \end{CD}$$
           the first and last term are vanishing because of the dimension assumption $b> b_1+1$.
           So the map $\mu $ can be described by a stable map of pairs. Now apply equivariant
           cohomology to this map of pairs.
            Since the ${\T}$-action on $(TF,TF^{{\T}})$
            is relatively free, the equivariant cohomology group
            $ H^*_{{\T}}(TF,TF^{{\T}}) $ identifies with
            the nonequivariant cohomology $H^*(TF/{{\T}},TF^{{\T}})$ of the quotient, which after replacing
            $TF^{{\T}}$ by a tubular neighbourhood, is a connected
            manifold relative to its boundary. 
             An orientation 
            of $H^1(X;{\R})$ and thus of $Pic^0(X)$ together with the standard orientation 
            of complex vector bundles
            defines an orientation class $[TF/{{\T}}]$ of this manifold.   
            Similarly, 
            the chosen homology orientation of $X$ and the 
            orientation of $Pic^0(X)$ determine 
            the orientation of $V$ and thus a generator $[V^+]$ in    
            reduced equivariant cohomology of $V^+$ as a free 
            $H^*_{{\T}}(*)\cong{\Z}[x]$-module of rank one. 
            The homomorphism $t$ associates
            to $\mu$ the degree zero part of 
            $\mu^*(\sum_0^\infty x^i[V^+])\cap [TF/{{\T}}]$. Using the alternate 
            description above, the same integer is obtained as
             $(\gamma({\C})^k\circ\mu)^*([(V\oplus{\C}^k)^+])\cap [TF/{{\T}}]$.

\end{pf}

\noindent
If the first Betti number of $X$ vanishes, the group 
$\pi_{{\T},{\CC}}^b(Pic^0(X); ind({\cal D}))$ simplifies:
The index of the Dirac operator is a complex vector space of complex 
dimension \[d={c({\fs})^2-signature(X)\over 8},\] where $c({\fs})$ is the first
Chern class of the spinor bundles $S^\pm$ associated 
to the $\sc$-structure $\fs$.

\begin{prop} For $i>1$,  the 
             stable equivariant cohomotopy groups $\pi_{{\T},{H}}^i(*; {\C}^d))$
             are isomorphic to the nonequivariant stable cohomotopy groups
             $\pi^{i-1}({\C}P^{d-1})$ of complex projective $(d-1)$-space.
             In particular, if $X$ is a closed 4-manifold with $b_1=0$ and $b_+>1$,
             then the monopole map determines an element in 
             $\pi^{b-1}({\C}P^{d-1})$.
\end{prop}            

\begin{pf}   The long exact stable cohomotopy sequence for the  pair $(D({\C}^d), S({\C}^d))$
             consisting of the unit disk and sphere in the complex vector space ${\C}^d$
             allows to identify for $i>1$ the groups $\pi_{{\T},H}^i(*; {\C}^d))$
             with $\tilde\pi_{{\T},H}^{i-1}(S({\C}^d)^+)$. But for the free ${\T}$-space 
             $S({\C}^d)$ equivariant cohomotopy is isomorphic to the nonequivariant 
             cohomotopy of its quotient \cite{tD}. 
\end{pf}

\noindent 
To analyse this cohomotopy of projective spaces a little further, consider the Hurewicz map 
\begin{eqnarray*}
       \pi^i(Y)&\to& H^i(Y)\cr
       [f]&\mapsto&f^*(1),
\end{eqnarray*}
with $1\in H^i(S^i)\cong \tilde H^0(S^0)$ defined by the orientation. 
Rationally it is an isomorphism,
as rationally the sphere spectrum is an Eilenberg MacLane spectrum by Serre's theorem. However, 
nonrationally, the Hurewicz map has both kernel and cokernel, as displayed in low dimensions
below.

In the following lemma the results are ordered according to k, which can be interpreted
as the "expected dimension of the moduli space", i.e. the dimension of the preimage of a generic
point in the sphere. Thanks go to N. Minami for pointing out a mistake in the computation for $k=3$
in an earlier version.

\begin{lem}\label{Hurewicz} Let $d>1$ be an integer. The Hurewicz map of reduced cohomology groups
            $h^{2d-2-k}:\widetilde{\pi}^{2d-2-k}({\C}P^{d-1})\to \widetilde{H}^{2d-2-k}({\C}P^{d-1})$
           \begin{enumerate}\setcounter{enumi}{-1}
         \item  for $k=0$ is an isomorphism.
         \item  for $k=1$ has kernel isomorphic to ${\Z}/{gcd(2,d)}$.
         \item  for $k=2$ has kernel isomorphic to ${\Z}/{gcd(2,d)}$ and cokernel to 
                ${\Z}/{gcd(2,d-1)}$.
         \item  for $k=3$ has kernel isomorphic to ${\Z}/l$ with $l= gcd(24,d)$, if $d$ is even, and 
                $l=gcd(24,d-3)/2$ otherwise.
         \item  for $k=4$ has trivial kernel and, for $d>2$, cokernel isomorphic to ${\Z}/m$ with
                $lm=48$, if $d$ is even, and $lm=12$ otherwise.
           \end{enumerate}
\end{lem}
   
\begin{pf} The proof employs the Atiyah-Hirzebruch spectral sequence with
           $E_2$-term 
           \[
              H^*({Y};\pi^*(pt))\Rightarrow \pi^*(Y)
           \] 
           and uses the following facts:\hfill\break
           1. The attaching map of the 4-cell in ${\C}P^2$ is 
           the Hopf map, which is the generator $\eta$ of $\pi^{-1}(pt)\cong {\Z}/2$.
           \hfill\break
           2. The group $\pi^{-2}(pt)\cong {\Z}/2$ is generated by $\eta^2$.\hfill\break
           3. The attaching map of the 8-cell in ${\H}P^2$ is again a Hopf map, which  is
           stably the generator $\nu$ of $\pi^{-3}(pt)\cong {\Z}/24$. Furthermore,
           $\eta^3=12\nu$.\hfill\break
           4. The stable homotopy groups $\pi^{-4}(pt)$ and $\pi^{-5}(pt)$ vanish.
           \hfill\break
           5. For even $d$, there is a projection of the complex projective to  the 
           quaternionic projective space.\hfill\break
           6. The differentials in the spectral sequence are differentials for the
           algebra structure on the respective $E_i$-terms.\hfill\break
           7. The spectral sequence is natural in $Y$. In particular 
           one may use the map between spectral sequences induced by
           inclusions of the projective spaces into higher dimensional
           projective spaces and induced by 
           the projection of complex projective 
           spaces to quaternionic projective spaces. \hfill\break
           8. By a result of I. M. James, cf. \cite{James}, the projection map 
           \[{\C}P^{d-1}/{\C}P^{d-4}\to {\C}P^{d-1}/{\C}P^{d-2}=S^{2d-2}\]
           is stably split if and only if $d$ is divisible by $24$. \hfill\break
\end{pf}

Let $m(d,k)$ denote the order of the cokernel of the Hurewicz map
$h^{2d-2-k}$ for even integers $k$. Since the integer valued Seiberg-Witten invariants are in the
image of the Hurewicz map, one gets on particular:

\begin{cor} For a $K$-oriented 4-manifold $X$ with vanishing first Betti number and $b^+=2p+1$ odd, the 
integer valued Seiberg-Witten invariant is divisible by $m(d,k)$ with $k=2d-2p-2$.
\end{cor}

Indeed, one gets estimates for the divisibility $m(d,k)$ of Seiberg-Witten invariants by comparing
the Hurewicz maps from stable cohomotopy to $K$-theory and to singular cohomology. The following 
result was obtained in \cite{Furuta2}. 

\begin{thm}\label{divisibility} The integers $m(d,2\kappa)$ are 
divisible by the denominators of the rational numbers
$a_{p,0}, a_{p,1}, \dots, a_{p, \kappa}$, which appear as coefficients 
in the Taylor expansion
\[log(1-\xi)^p=\sum_{l=0}^{\infty}a_{p,l}\xi^{p+l}.\]
\end{thm}

Comparison with \ref{Hurewicz} shows that the above estimates are not sharp in general.

\begin{pf} For a stable map $\mu:{\C}P^{d-1}\to S^{2p}$, consider the commuting
diagram
\[\begin{CD}  K(S^{2p}) & @>{\mu^*}>> & K({\C}P^{d-1})\\
               @V{ch}V{\cong}V &              & @VV{ch}V           \\
              H^*(S^{2p};{\Z}) & @>{\mu^*}>> & H^*({\C}P^{d-1};{\Q})
\end{CD}  
\]
with vertical arrows given by the Chern character and horizontal maps induced by $\mu$. 
The left vertical map is an isomorphism,
mapping the class $\tau$ to the orientation class $[S^{2p}]\in H^*(S^{2p};{\Z})$.
The Hurewicz image in 
$ K({\C}P^{d-1})\cong {\Z}[\xi]/(\xi^d)$ of the cohomotopy element $\mu$ is 
of the form $\mu^*(\tau)=\sum_{l=0}^{\kappa}b_l \xi^{l +p}$ with integers $b_l$, because the composite
of $\mu$ with the inclusion map ${\C}P^{p-1}\to {\C}P^{d-1} $ is nullhomotopic by dimension reasons.
The Hurewicz image of $\mu$ in $ H^*({\C}P^{d-1};{\Q})\cong {\Q}[x]/(x^d)$ is of the form
$nx^p$ for an integer $n$.
The Chern character $\xi\mapsto (1-exp(x))$ on the right hand side is injective and becomes an isomorphism 
after tensoring with the rationals. Commutativity of the diagram now implies
\[
\sum_{l=0}^{\kappa}b_l \xi^{l +p}= \mu^*(\tau) = ch^{-1}\mu^*([S^{2p}])=ch^{-1}(nx^p)=d\, log(1-\xi)^p
=n\sum_{l=0}^{\kappa}a_{p,l}\xi^{p+l}.
\]
\end{pf}

To show that the 
stable cohomotopy invariants 
are indeed effective refinements of the integer valued Seiberg-Witten invariants,
it remains to find manifolds whose stable cohomotopy invariants
are elements in the kernel of the Hurewicz-map. This will be done
in \cite{PartII}.

\noindent 
Let's finish by showing how to recapture Donaldon's first theorem
on manifolds with definite intersection form in the stable cohomotopy
setting. This proof relies on
the following fact about equivariant maps, which is well known and 
can be proved
basically the same way as \ref{divisibility}
by the use of the equivariant $K$-theory mapping degree (cf. e.g. \cite{tD}):

\begin{lem}\label{Tammo}
Let $f:({\R}^n\oplus {\C}^m)^+\to ({\R}^n\oplus {\C}^{m+k})^+$ be an ${\T}$-equivariant
map such that the restricted map on the fixed points has degree $1$. Then
$k\geq 0$ and $f$ is homotopic to the inclusion. 
\end{lem}

\noindent
We apply this lemma to the case where $b=b_+=0$, i.e. to manifolds
$N$ with negative definite intersection form. Note that for these 
manifolds the 
monopole map on the fixed point set is just a linear isomorphism in each fiber
over $Pic^0(N)$.
As a consequence one gets:

\begin{cor} (Donaldson\cite{Don}) 
Let $N$ be a closed oriented four-manifold with negative 
definite intersection form. Then every characteristic $c\in H_2(N;{\Z})$
satisfies $-c^2\geq b_2(N)$. As a consequence, by \cite{Elkies} the 
intersection pairing is diagonal.
\end{cor}

\begin{pf} The stable cohomotopy invariant associated to a $\sc$-structure
with determinant $c$ is in $\pi^0_{{\T},H}(Pic^0(N); {\mathrm ind}(D))$.
Restricting it to a point in $Pic^0(N)$ results in 
the stabilisation of a map as in \ref{Tammo}
with \[k=-{\mathrm ind}_{\C}(D)= \frac{-c^2-b_2(N)}{8}.\]
\end{pf}


\end{document}